\documentclass[12pt,reqno]{article}

\usepackage[usenames]{color}
\usepackage[colorlinks=true,
linkcolor=webgreen, filecolor=webbrown,
citecolor=webgreen]{hyperref}

\definecolor{webgreen}{rgb}{0,.5,0}
\definecolor{webbrown}{rgb}{.6,0,0}

\usepackage{amssymb}
\usepackage{graphicx}
\usepackage{amscd}
\usepackage{lscape}
\usepackage{tikz}
\usepackage{tikz-cd}
\usepackage{pgfplots}

\usetikzlibrary{matrix}
\usetikzlibrary{fit,shapes}
\usetikzlibrary{positioning, calc}
\tikzset{circle node/.style = {circle,inner sep=1pt,draw, fill=white},
        X node/.style = {fill=white, inner sep=1pt},
        dot node/.style = {circle, draw, inner sep=5pt}
        }
\usepackage{tkz-fct}

\usepackage{amsthm}
\newtheorem{theorem}{Theorem}

\newtheorem{conjecture}[theorem]{Conjecture}

\theoremstyle{definition}

\newtheorem{example}[theorem]{Example}

\usepackage{float}

\usepackage{graphics,amsmath}
\usepackage{amsfonts}
\usepackage{latexsym}
\usepackage{epsf}

\newcommand{\seqnum}[1]{\href{http://oeis.org/#1}{\underline{#1}}}

\begin{document}

\begin{center}
\vskip 1cm{\LARGE\bf Conjectures on Somos $4$, $6$ and $8$ sequences using Riordan arrays and the Catalan numbers} \vskip 1cm \large
Paul Barry\\
School of Science\\
South East Technological University\\
Ireland\\
\href{mailto:pbarry@wit.ie}{\tt pbarry@wit.ie}
\end{center}
\vskip .2 in

\begin{abstract} We give conjectures on the form of families of integer sequences whose Hankel transforms are, respectively, $(\alpha, \beta)$ Somos $4$ sequences, $(\alpha, 0, \gamma)$ Somos $6$ sequences, and $(\alpha, \beta, \gamma, \delta)$ Somos $8$ sequences, for particular values of $\alpha$, $\beta$, $\gamma$, $\delta$ which we describe. The sequences involved can be described in terms of the application of certain stretched Riordan arrays to the Catalan numbers, accompanied by a (sequence) Hankel transform. The combination of Riordan array and the Catalan numbers results from the study of certain generalized Jacobi continued fractions, based on the Counting Automata Methodology.
\end{abstract}

There is a growing literature on Somos sequences. Most notably, Hone, and collaborators, in a series of papers \cite{Braden, Fed, Hamad, Shadow, Hone1, Hone2, Hone3, Hone4, Hone5}, has given a comprehensive perspective on Somos sequences. This follows work by van Poorten \cite{Poorten1, Poorten2}. Hone exhibits Somos sequences in the context of special integrable systems, where the tools of Lax operators and tau functions are particularly relevant. In their specific details, these sequences can be defined in terms of elliptic and hyperelliptic sigma functions.

As part of the arsenal of tools that can be employed to derive Somos sequences, the (sequence) Hankel transform has been employed \cite{Chang1, Chang2, Shadow, Hone1, Hankel, Xin}. In particular, the theory of Jacobi continued fractions \cite{Wall} and orthogonal polynomials has proved useful.

In this note, we indicate how a type of generalized Jacobi continued fraction can lead to conjectures concerning the generation of Somos $4$ and Somos $6$ sequences as the Hankel transforms of integer sequences, with these latter sequences being expressible in terms of the Catalan numbers and certain (stretched) Riordan arrays \cite{book, bbook, SGWW}. These generalized Jacobi continued fractions appear notably in work by de Castro, Ram\'irez, and Ram\'irez using their Counting Automata Methodology \cite{DRR}. The following example shows how these generalized Jacobi continued fractions can be applied to lattice paths whose steps are weighted according to a weight function.
\begin{example} \label{ex1} The sequence $g_n$, with generating function given by the generalized Jacobi continued fraction
$$g(x)=\cfrac{1}{1-\frac{x}{1-x}-\cfrac{x^2}{1-\frac{x}{1-x}-\cfrac{x^2}{1-\frac{x}{1-x}-\cdots}}},$$ which begins
$$1, 1, 3, 7, 19, 51, 143, 407, 1183, 3487, 10415,\ldots,$$ counts the number of lattice paths, never going below the $x$-axis, from $(0,0)$ to $(n,0)$, consisting of up steps $U = (1,1)$, down steps $D = (1,-1)$ and horizontal steps $H(k) = (k,0)$ for every positive integer $k$.
To obtain a closed form for $g(x)$, we solve the equation
$$u=\frac{1}{1-\frac{x}{1-x}-x^2u}$$ for $u$ (with $u(0)=1$), to get
$$g(x)=\frac{1-2x-\sqrt{1-4x+8x^3-4x^4}}{2x^2(1-x)}.$$ Using the generating function $c(x)=\frac{1-\sqrt{1-4x}}{2x}$ of the Catalan numbers $C_n=\frac{1}{n+1}\binom{2n}{n}$ \seqnum{A000108}, we can express this as
$$g(x)=\frac{1-x}{1-2x}c\left(\frac{x^2(1-x)^2}{(1-2x)^2}\right).$$
This means that the expansion $1,1,3,7,19,51,143,\ldots$ \seqnum{A135052} of $g(x)$ can be obtained by applying the (stretched) Riordan array $\left(\frac{1-x}{1-2x},\frac{x^2(1-x)^2}{(1-2x)^2}\right)$ to (the vector of) the Catalan numbers:
$$\left(
\begin{array}{cccccccc}
 1 & 0 & 0 & 0 & 0 & 0 & 0 & 0 \\
 1 & 0 & 0 & 0 & 0 & 0 & 0 & 0 \\
 2 & 1 & 0 & 0 & 0 & 0 & 0 & 0 \\
 4 & 3 & 0 & 0 & 0 & 0 & 0 & 0 \\
 8 & 9 & 1 & 0 & 0 & 0 & 0 & 0 \\
 16 & 25 & 5 & 0 & 0 & 0 & 0 & 0 \\
 32 & 66 & 20 & 1 & 0 & 0 & 0 & 0 \\
 64 & 168 & 70 & 7 & 0 & 0 & 0 & 0 \\
\end{array}
\right)\left(\begin{array}{c} 1\\1\\2\\5\\14\\42\\132\\429\\\end{array}\right)=
\left(\begin{array}{c} 1\\1\\3\\7\\19\\51\\143\\407\\\end{array}\right).$$
We have $g_n=\sum_{k=0}^n \left(\sum_{j=0}^{n-2k} \binom{2k+1}{j}(-1)^j \binom{n-j}{2k}2^{n-2k-j}\right)C_k$. The Hankel transform $h_n=|g_{i+j}|_{0 \le i,j \le n}$ is then given by
$$h_n = 2^{\lfloor \frac{(n+1)^2}{4} \rfloor}.$$
This Hankel transform $h_n$ is then a $(0,4)$ Somos $4$ sequence, since we have, for $n>4$,
$$h_n = 4 \frac{h_{n-2}^2}{h_{n-4}}.$$

More generally, if 
$$g(x)=\frac{1}{1-\frac{x}{1-rx}-x^2 g(x)},$$ 
then we have 
$$g(x)=\frac{1-rx}{1-(r+1)x}c\left(\frac{x^2(1-rx)^2}{(1-(r+1)x)^2}\right).$$ This expands to give the sequence $g_n$ with 
$$g_n=\sum_{k=0}^n \left(\sum_{j=0}^{n-2k} \binom{2k+1}{j}(-r)^j \binom{n-j}{2k}(r+1)^{n-2k-j}\right)C_k.$$
The Hankel transform of $g_n$ is conjectured to be a $((r-1)^2, 4r)$ Somos $4$ sequence. Note that the sequence 
$$\tilde{g}_n=\sum_{k=0}^n \left(\sum_{j=0}^{n-2k} \binom{2k+1}{j}(-1)^j \binom{n-j}{2k}(r+1)^{n-2k-j}\right)C_k$$
will have the same Hankel transform. 
\end{example}

\section{Somos sequences}
In this note, we deal with Somos $4$ and Somos $6$ sequences. An $(\alpha, \beta)$ Somos $4$ sequence $a_n$ is a sequence such that 
$$a_n = \frac{\alpha a_{n-1}a_{n-3} + \beta a_{n-2}^2}{a_{n-4}}, n>4,$$ for appropriate initial values.

An $(\alpha, \beta, \gamma)$ Somos $6$ sequence $a_n$ is a sequence such that 
$$a_n = \frac{\alpha a_{n-1}a_{n-5} + \beta a_{n-2} a_{n-4} + \gamma a_{n-3}^2}{a_{n-6}}, n>6,$$ for appropriate initial values.

An $(\alpha, \beta, \gamma, \delta)$ Somos $8$ sequence $a_n$ is a sequence such that
$$a_n = \frac{\alpha a_{n-1}a_{n-7} + \beta a_{n-2} a_{n-6} +\gamma a_{n-3} a_{n-5}+ \delta a_{n-4}^2}{a_{n-8}}, n>8,$$ for appropriate initial values.

\section{Riordan arrays} A Riordan array $(g(x), f(x))$ may be defined by two power series 
$$g(x)=g_0 + g_1 x + g_2 x^2+ \cdots, g_0 \ne 0,$$ and 
$$f(x)=f_1 x+ f_2 x^2+ f_3x^3 + \cdots, f_0=0, f_1 \ne 0,$$ with the coefficients $g_n$ and $f_n$ drawn from an appropriate field (or from an appropriate ring: $\mathbb{Z}$ is the relevant ring for this note). The term ``array'' comes from the matrix representation of a Riordan array, which is the matrix $(t_{n,k})_{0 \le n,k \le \infty}$ where
$$t_{n,k}=[x^n] g(x)f(x)^k.$$ Here, $[x^n]$ is the functional that extracts the coefficient of $x^n$ in a power series. A Riordan array is a lower triangular invertible matrix. 

A  Riordan array $(g(x), f(x))$ acts on a power series $h(x)$ by the action (weighted composition)
$$(g(x), f(x))h(x)= g(x)h(f(x)).$$ 
In matrix terms, this is equivalent to multiplying the vector $(h_n)$, where $h(x)=\sum_{n=0}^{\infty} h_n x^n$, by the matrix $(t_{n,k})$. 

The matrix representation of couples of the form $(g(x), x^r f(x))$ where $f(x)$ is as above, and the integer $r>0$, gives what are called "stretched" Riordan arrays. An example has already been met in Example (\ref{ex1}).

\section{Jacobi continued fractions and Hankel transforms}
A continued fraction of the form 
$$g(x)=\cfrac{1}{1- a_0 x- \cfrac{b_0 x^2}{1-a_1 x - \cfrac{b_1 x^2}{1- a_2 x-\cdots}}}$$ is called a Jacobi continued fraction. We note that the Hankel transform $n_n=|g_{i+j}|_{0 \le i,j \le \infty}$ of the expansion of such a continued fraction will be given by $$h_n = \prod_{k=0}^n b_k^{n-k}.$$ This is independent of the coefficients $a_n$. 

\section{Generalized Jacobi continued fractions and Somos $4$ sequences}
In this section, we itemize in increasing detail conjectures concerning generalized Jacobi continued fractions and Somos $4$ sequences. The examples given are in increasing order of complexity. The goal is to produce integer sequences whose Hankel transform are $(\alpha, \beta)$ Somos $4$ sequences.
\begin{conjecture} We consider the generalized Jacobi continued fraction 
$$g(x)=\frac{1}{1-\frac{1+rx}{1-x} x - sx^2 g(x)}.$$ Then we can express $g(x)$ as 
$$g(x)=\frac{1-x}{1-2x-rx^2}c\left(\frac{sx^2(1-x)^2}{(1-2x-rx^2)^2}\right).$$ 
The Hankel transform $h_n$ of the sequence $g_n$ is given by
$$h_n = s^{\lfloor \frac{n^2}{4} \rfloor}(r+s+1)^{\lfloor \frac{(n+1)^2}{4} \rfloor}.$$ 
Then $h_n$ is a $(0, s^2(r+s+1)^2)$ Somos $4$ sequence.
\end{conjecture}
The generating function $g(x)$ is obtained by applying the (stretched) Riordan array 
$$\left(\frac{1-x}{1-2x-rx^2}, \frac{sx^2(1-x)^2}{(1-2x-rx^2)^2}\right)$$ to the generating function $c(x)$ of the Catalan numbers. 
\begin{conjecture}  We consider the generalized Jacobi continued fraction
$$g(x)=\frac{1}{1-\frac{1+rx}{1-x} x - \frac{sx^2}{1-x} g(x)}.$$ We can express $g(x)$ as 
$$g(x)=\frac{1-x}{1-2x-rx^2}c\left(\frac{sx^2(1-x)}{(1-2x-rx^2)^2}\right).$$ 
The Hankel transform $h_n$ of the sequence $g_n$ is a $(s^2, s^2(r+(r+s)^2))$ Somos $4$ sequence.
\end{conjecture}
\begin{conjecture} We consider the generalized Jacobi continued fraction
$$g(x)=\frac{1}{1-\frac{1+rx}{1-x} x - \frac{1+sx}{1-x} x^2g(x)}.$$ We can express $g(x)$ as
$$g(x)=\frac{1-x}{1-2x-rx^2}c\left(\frac{x^2(1-x)(1+sx)}{(1-2x-rx^2)^2}\right).$$
The Hankel transform $h_n$ of the sequence $g_n$ is a $((s+1)^2,(1+r^2-6s-3s^2-r(s^2+2s-3)))$ Somos $4$ sequence.
\end{conjecture}
The last two conjectures are encompassed in the following more general conjecture.
 \begin{conjecture} We consider the generalized Jacobi continued fraction
$$g(x)=\frac{1}{1-v\frac{1+rx}{1-x} x - w\frac{1+sx}{1-x} x^2g(x)}.$$ Then we can express $g(x)$ as
$$\frac{1-x}{1-(v+1)x-rx^2}c\left(\frac{wx^2(1-x)(1+sx)}{(1-(v+1)x-rx^2)^2}\right).$$ 
The Hankel transform $h_n$ of $g_n$ is an $(\alpha, \beta)$ Somos $4$ sequence with parameters 
$$\alpha = (s+v)^2 w^2,$$ and 
$$\beta=w^2(r^2v^2+w(w+v-v^2)+rv(v+2w)-s^2(v(r+1)+2w)-s((r+1)v^2+w+v(r+1+3w))).$$ 
\end{conjecture}

\section{Generalized Jacobi continued fractions and Somos $6$ sequences}
We now extend the ideas of the last section to formulate some conjectures concerning the Hankel transform of integer sequences and Somos $6$ sequences of type $(\alpha,0,\gamma)$. We start with an illustrative example.
\begin{example} We consider the generalized Jacobi continued fraction 
$$g(x)=\frac{1}{1-x \frac{1+3x}{1-x}+ x^2 \frac{1+2x}{1-x}-x^3 g(x)},$$ which has closed form 
$$g(x)=\frac{1-2x-2x^2+2x^3-\sqrt{1-4x+8x^3+4x^4-12x^5+4x^6}}{2x^3(1-x)}$$ or in Catalan form, as
$$\frac{1-x}{1-2x-2x^2+2x^3}c\left(\frac{x^3(1-x)^2}{(1-2x-2x^3+2x^3)^2}\right).$$ 
The generating function $g(x)$ expands to give the sequence $g_n$ that begins
$$1,1,4,9,25,67,183,\ldots$$ whose Hankel transform $h_n$ begins 
$$1,3,2,-23,-231,-1987,-41482,\ldots.$$ We then conjecture that $h_n$ is a $(9,0,23)$ Somos $6$ sequence, that is, we have 
$$h_n = \frac{9 e_{n-1} e_{n-5}+23 e_{n-3}^2}{e_{n-6}}, \quad n>6.$$ 
\end{example}
\begin{conjecture} We consider the generalized Jacobi continued fraction 
$$g(x)=\frac{1}{1-x\frac{1+rx}{1-x}-x^2\frac{1+sx}{1-x}-tx^3 g(x)}.$$ We have that 
$$g(x)=\frac{1-2x-(r+1)x^2-sx^3-\sqrt{Q(x,r,s,t)}}{2tx^3(1-x)},$$ where 
$$Q(x,r,s,t)=1-4x-2(r-1)x^2+2(2r-s-2(t-1))x^3+(r^2+2r+4s+8t+1)x^4$$
$$\quad\quad\quad +2(rs+s-2t)x^5+s^2x^6.$$ 
In Catalan form, we have 
$$g(x)=\frac{1-x}{1-2x-(r+1)x^2-sx^3}c\left(\frac{tx^3(1-x)^2}{(1-2x-(r+1)x^2-sx^3)^2}\right).$$ 
We conjecture that the Hankel transform $h_n$ of $g_n$ is an $(\alpha, 0,\gamma)$ Somos $6$ sequence with parameters
\begin{align*}
\alpha &=t^2(r+2)^2\\
\gamma &={\scriptstyle t^3(r^3t+r^2(s+7t)+2r(s^2+2(t+1)s+t(t+8))+s^3+s^2(3t+4)+s(t+2)(3t+2)+t(t^2+4t+12))}.\end{align*}
\end{conjecture}
The Riordan array 
$$\left(\frac{1-x}{1-2x-(r+1)x^2-sx^3}, \frac{tx(1-x)^2}{(1-2x-(r+1)x^2-sx^3)^2}\right)$$ has its general term $t_{n,k}=t_{n,k}(r,s,t)$ given by 
$$t^k \sum_{j=0}^{2k+1}\binom{2k+1}{j}(-1)^j \sum_{i=0}^{n-k-j}\binom{2k+i}{i}\sum_{m=0}^i 2^{i-m}\binom{m}{n-j-i-m} s^{n-k-j-i-m}(r+1)^{2m-n+k+j+i}.$$ 
Then we have 
$$g_n=\sum_{k=0}^{\lfloor \frac{n}{3} \rfloor} t_{n-2k,k}C_k.$$ 
Assuming the validity of the conjecture, this then gives us a mechanism for producing a three parameter family of integer sequences whose Hankel transforms are Somos $6$ sequences. 
\begin{example} The sequence $g_n(-2,-2,1)=\sum_{k=0}^{\lfloor \frac{n}{3} \rfloor} t_{n-2k,k}(-2,-2,1)C_k$ begins $$1, 1, 0, -3, -7, -9, -5, 8, 32, 71, 129, 187, 153, \ldots.$$ 
This has a Hankel transform that begins 
$$1, -1, -2, 5, 17, -3, -122, 1201, -2980,\ldots.$$ The conjecture is that this is a $(1,0,-5)$ Somos $6$ sequence.
\end{example}
\begin{example} Taking $r=-3, s=0, t=-1$ we obtain a sequence that begins 
$$1, 1, 0, -3, -7, -7, 7, 42, 78, 35, -217, -695, -907, 523, \ldots.$$ The Hankel transform of this sequence begins $$1, -1, -2, -3, 11, 23, 4, -355, -1326,\ldots.$$ The conjecture is that this is a $(1,0,3)$ Somos $6$ sequence.
\end{example}
We can add two more parameters by considering 
$$g_n(r,s,t,u,v)=\sum_{k=0}^{\lfloor \frac{n}{3} \rfloor} t_{n-2k,k}(r,s,t)u^k v^{n-2k}C_k.$$ 
Again, we conjecture that the Hankel transform of the sequence $g_n$ is a Somos $6$ sequence.
\begin{example}
We take $r=s=t=1$ and $u=2, v=-1$, to obtain the sequence $g_n$ that begins 
$$1, -1, 3, -10, 26, -75, 224, -659, 1979, -6025, 18452, -57028, 177625, \ldots,$$ with a Hankel transform $h_n$ that begins 
$$1, 2, -15, -182, -4864, 85976, 26865504, 5387832064, 687205582336,\ldots.$$
We conjecture that this is a $(16,0,728)$ Somos $6$ sequence. 
\end{example}

\section{Conjectures on integer Somos $8$ sequences}
The classical Somos $8$ sequence beginning $1,1,1,1,1,1,1,1,\ldots$ with parameters $(1,1,1,1)$ is not an integer sequence \cite{Laurent}. In fact, it begins 
 $$1,1,1,1,1,1,1,1,4,7,13,25,61,187,775,5827,14815,\frac{420514}{7},\frac{28670773}{91},\frac{6905822101}{2275},\ldots.$$ 
 In this section, using Hankel transforms of integer sequences, we conjecture the form of infinite families of integer Somos $8$ sequences. In general, the $(\alpha, \beta, \gamma, \delta)$ parameters are rational numbers. This is the content of the following conjectures. The integrality of the sequences arises as we consider the Hankel transforms of integer sequences. 
\begin{conjecture} We consider the continued fraction 
$$g(x)=\frac{1}{1-\frac{x}{1-rx}-x^2-x^3 g(x)}.$$ 
We have
$$g(x)=\frac{1-rx}{1-(r+1)x-x^2+rx^3}c\left(\frac{x^3(1-rx)^2}{(1-(r+1)x-x^2+rx^3)^2}\right).$$
We have \begin{scriptsize}
$$g_n=\sum_{k=0}^n (\sum_{j=0}^{2k+1} \binom{2k+1}{j}(-r)^j \sum_{i=0}^{n-3k} \binom{2k+i}{i}\sum_{m=0}^i \binom{i}{m}(r+1)^{i-m} \binom{m}{n-3k-j-i-m}(-r)^{n-3k-j-i-m})C_k.$$ 
\end{scriptsize}
The Hankel transform $h_n$ of the sequence $g_n$ is an integer $(\alpha, \beta, \gamma, \delta)$ Somos $8$ sequence with
\begin{align*}
\alpha&=-\frac{-r^8+8 r^7-21 r^6+40 r^5-35 r^4+24 r^3-71 r^2-8 r}{r^4-2 r^3+8 r^2+2 r-9}\\
\beta&=\frac{8 \left(r^9-6 r^8+17 r^7-30 r^6+15 r^5-14 r^4-r^3-14 r^2\right)}{r^4-2 r^3+8 r^2+2 r-9}\\
\gamma&=\frac{8 \left(r^{10}-2 r^8+29 r^7-32 r^6+39 r^5+18 r^4+11 r^3-r^2+r\right)}{r^3-r^2+7 r+9}\\
\delta&={\scriptstyle -\frac{-2 r^{13}+13 r^{12}-48 r^{11}+85 r^{10}-83 r^9+11 r^8-124 r^7+454 r^6-364 r^5+263 r^4+84 r^3+189 r^2+25 r+9}{r^4-2 r^3+8 r^2+2 r-9}}.\end{align*}
\end{conjecture}
\begin{conjecture}
We consider the continued fraction 
$$g(x)=\frac{1}{1-x-\frac{x^2}{1-rx}-x^3g(x)}.$$ 
We have 
$$g(x)=\frac{1-rx}{1-(r+1)x+(r-1)x^2}c\left(\frac{x^3(1-rx)^2}{(1-(r+1)x+(r-1)x^2)^2}\right).$$ 
Then \begin{scriptsize}
$$g_n=\sum_{k=0}^n (\sum_{j=0}^{2k+1}\binom{2k+1}{j}(-r)^j \sum_{i=0}^{n-3k} \binom{2k+i}{i}\binom{i}{n-3k-j-i}(1-r)^{n-3k-j-i}(r+1)^{2i-n+3k+j})C_k.$$
\end{scriptsize}
The Hankel transform $h_n$ of the sequence $g_n$ is an integer $(\alpha, \beta, \gamma, \delta)$ Somos $8$ sequence with
\begin{align*}
\alpha &=-\frac{-r^4+11 r^3-26 r^2+16 r+5}{r^2-4 r+3}\\
\beta &=-\frac{-2 r^5+19 r^4-40 r^3+13 r^2+5 r}{r^2-4 r+3}\\
\gamma &=-\frac{-3 r^6+12 r^5-15 r^4-25 r^3+62 r^2+36 r+5}{r-3}\\
\delta &=-\frac{-r^9+8 r^8-26 r^7+43 r^6-40 r^5+17 r^4+23 r^3-27 r^2-19 r-3}{r^2-4 r+3}.\end{align*}
\end{conjecture}
\begin{conjecture} We consider the continued fraction 
$$g(x)=\frac{1}{1-\frac{1-(r-1)x}{1-rx}x-x^2-x^3 g(x)}.$$ 
We have 
$$g(x)=\frac{1-rx}{1+(r+1)x+(r-2)x^2+rx^3}c\left(\frac{x^3 (1-rx)^2}{(1+(r+1)x+(r-2)x^2+rx^3)^2}\right).$$ 
Then \begin{scriptsize}
$$g_n=\sum_{k=0}^n (\sum_{j=0}^{2k+1} \binom{2k+1}{j}(-r)^j \sum_{i=0}^{n-3k} \sum_{m=0}^i \binom{i}{m}(-1)^m (r+1)^{i-m} \binom{m}{n-3k-j-i-m}r^{n-3k-j-i-m}(r-2)^{2m-n+3k+j+i})C_k.$$ 
\end{scriptsize}
The Hankel transform $h_n$ of the sequence $g_n$ is an integer $(\alpha, \beta, \gamma, \delta)$ Somos $8$ sequence with 
\begin{align*}
\alpha &=-\frac{r^7-8 r^6+25 r^5-20 r^4-37 r^3+75 r+28}{2 \left(r^3-3 r^2-5 r+7\right)}\\
\beta &=\frac{(r+1) \left(r^8-11 r^7+47 r^6-83 r^5+17 r^4+71 r^3+45 r^2-169 r+210\right)}{2 \left(r^3-3 r^2-5 r+7\right)}\\
\gamma &=\frac{(r^2-1) \left(3 r^8-29 r^7+115 r^6-225 r^5+181 r^4+105 r^3-255 r^2-235 r+84\right)}{2 \left(r^3-3r^2-5r+7\right)}\\
\delta &=-\frac{r^{10}-17 r^9+96 r^8-212 r^7+54 r^6+594 r^5-796 r^4-36 r^3+721 r^2-329 r-588}{2 \left(r^3-3 r^2-5 r+7\right)}.\end{align*}
\end{conjecture}

\begin{example} The conjectures above are not exhaustive. We consider, for instance, the continued fraction 
$$g(x)=\frac{1}{1-\frac{x}{1-\frac{x}{1-3x}}-x^3 g(x)}=\frac{1}{1-\frac{x(1-3x)}{1-4x}-x^3 g(x)}.$$ 
We find that 
$$g(x)=\frac{1-5x+3x^2-\sqrt{1-10x+31x^2-34x^3+41x^4-64x^5}}{2x^3(1-4x)},$$ or equivalently,
$$g(x)=\frac{1-4x}{1-5x+3x^2}c\left(\frac{x^3(1-4x)^2}{(1-5x+3x^2)^2}\right).$$ 
This expands to give the integer sequence $g_n$ that begins
$$1, 1, 2, 8, 32, 133, 569, 2450, 10569, 45643, 197206, 852239, 3683553,\ldots.$$ 
This has a Hankel transform $h_n$ that begins 
$$1, 1, -8, -161, -1333, 631, 1570896, 194685449, 8871803329, -1552662557863, \ldots.$$ 
We now conjecture that the integer sequence $h_n$ is a $\left(-\frac{101}{3}, -\frac{484}{3}, 4299, \frac{23359}{3}\right)$ Somos $8$ sequence.
\end{example}
\begin{example} We consider the generating function $g(x)$ defined by 
$$g(x)=\frac{1}{1-x-\frac{x^2}{1-\frac{x}{1-x}}-x^3 g(x)}=\frac{1}{1-x-\frac{x^2(1-x)}{1-2x}-x^3 g(x)}.$$ 
We have 
$$g(x)=\frac{1-3x+x^2+x^3-\sqrt{1-6x+11x^2-8x^3+11x^4-14x^5+x^6}}{2x^3(1-2x)},$$ or equivalently,
$$g(x)=\frac{1-2x}{1-3x+x^2+x^3}c\left(\frac{x^3(1-2x)^2}{(1-3x+x^2+x^3)^2}\right).$$ 
The sequence $g_n$ begins 
$$1, 1, 2, 5, 12, 30, 77, 199, 518, 1357, 3572, 9443, 25064,\ldots,$$ with a Hankel transform that 
begins 
$$1, 1, -1, -4, -8, -13, 57, 241, 1093, 792, -30661, -246182,\ldots.$$ 
We conjecture that this integer sequence is a $\left(\frac{1}{2}, -\frac{5}{2}, \frac{11}{2}, \frac{17}{2}\right)$ Somos $8$ sequence.
\end{example}

\section{Polynomial sequences} 
We note that the parameterized sequences in the last section, whose Hankel transforms are conjectured to be Somos $8$ sequences, are in fact sequences of polynomials.
\begin{example} We consider the polynomial sequence (in $r$) with generating function 
$$\frac{1-rx}{1-(r+1)x+(r-1)x^2}c\left(\frac{x^3(1-rx)^2}{(1-(r+1)x+(r-1)x^2)^2}\right).$$ 
This expands to give the polynomial sequence that begins
$$1, 1, 2, 4 + r, 8 + 2 r + r^2, 17 + 5 r + 2 r^2 + r^3, 37 + 13 r +   6 r^2 + 2 r^3 + r^4,$$
 $$\quad\quad\quad 82 + 32 r + 16 r^2 + 7 r^3 + 2 r^4 + r^5,ldots.$$ 
The coefficient array of this polynomial sequence begins 
$$\left(
\begin{array}{ccccccccccc}
 1 & 0 & 0 & 0 & 0 & 0 & 0 & 0 & 0 & 0 & 0 \\
 1 & 0 & 0 & 0 & 0 & 0 & 0 & 0 & 0 & 0 & 0 \\
 2 & 0 & 0 & 0 & 0 & 0 & 0 & 0 & 0 & 0 & 0 \\
 4 & 1 & 0 & 0 & 0 & 0 & 0 & 0 & 0 & 0 & 0 \\
 8 & 2 & 1 & 0 & 0 & 0 & 0 & 0 & 0 & 0 & 0 \\
 17 & 5 & 2 & 1 & 0 & 0 & 0 & 0 & 0 & 0 & 0 \\
 37 & 13 & 6 & 2 & 1 & 0 & 0 & 0 & 0 & 0 & 0 \\
 82 & 32 & 16 & 7 & 2 & 1 & 0 & 0 & 0 & 0 & 0 \\
 185 & 80 & 41 & 19 & 8 & 2 & 1 & 0 & 0 & 0 & 0 \\
 423 & 201 & 108 & 51 & 22 & 9 & 2 & 1 & 0 & 0 & 0 \\
 978 & 505 & 282 & 140 & 62 & 25 & 10 & 2 & 1 & 0 & 0 \\
\end{array}
\right).$$ 
We note that the initial column of this array,
$$1, 1, 2, 4, 8, 17, 37, 82, 185, 423, 97,\ldots,$$
is the RNA sequence \seqnum{A004148}. The Hankel transform of the above polynomial sequence than begins 
$$1, 1, -2 r, -1 - 4 r - r^2 + 2 r^3 - r^4, -1 - 5 r - 6 r^2 + r^3 -   5 r^4 + 4 r^5 - r^6,$$
 $$\quad \quad \quad \quad -1 - 6 r - 7 r^2 + 12 r^3 + 12 r^4 + r^5 -   5 r^6 + r^7,\ldots.$$
\end{example}

\section{Conclusion} In this paper, we have conjectured that a combination of generating functions defined by generalized Jacobi continued fractions, Riordan arrays, Catalan numbers, and the sequence Hankel transform can be a fruitful context within which to explore Somos $4$, $6$ and $8$ sequences. It is noteworthy that if the conjectures about Somos $8$ sequences are true, then we can produce infinite families of integer Somos $8$ sequences.

\bigskip
\hrule
\bigskip
\noindent 2020 {\it Mathematics Subject Classification}: 
Primary 11B37; Secondary 05A15, 15B36, 11A55, 15A15, 15B05.
\noindent \emph{Keywords:} Somos sequence, generating function, Hankel transform, continued fraction, Riordan array, Catalan number.

\bigskip
\hrule
\bigskip
\noindent (Concerned with sequences
\seqnum{A000108}, 
\seqnum{A004148}, and
\seqnum{A135052}).

\end{document}